\documentstyle[12pt]{article}
\begin{document}
\bibliographystyle{tom}

\newtheorem{lemma}{Lemma}[section]
\newtheorem{thm}[lemma]{Theorem}
\newtheorem{cor}[lemma]{Corollary}
\newtheorem{voorb}[lemma]{Example}
\newtheorem{rem}[lemma]{Remark}
\newtheorem{prop}[lemma]{Proposition} 

\newfont{\gothic}{eufm10 scaled\magstep1}

\newcommand{\gota}{\mbox{\gothic a}}
\newcommand{\gotg}{\mbox{\gothic g}}
\newcommand{\gothh}{\mbox{\gothic h}}
\newcommand{\gotk}{\mbox{\gothic k}}
\newcommand{\gotm}{\mbox{\gothic m}}
\newcommand{\gotn}{\mbox{\gothic n}}
\newcommand{\gotp}{\mbox{\gothic p}}
\newcommand{\gotq}{\mbox{\gothic q}}
\newcommand{\gots}{\mbox{\gothic s}}

\newcommand{\sgotg}{{\mbox{\gothic g}}}

\newcounter{teller}
\renewcommand{\theteller}{\Roman{teller}}
\newenvironment{tabel}{\begin{list}%
{\rm \bf \Roman{teller}.\hfill}{\usecounter{teller} \leftmargin=1.1cm
\labelwidth=1.1cm \labelsep=0cm \parsep=0cm}
                      }{\end{list}}

\newcounter{proofstep}
\newcommand{\nextstep}{\refstepcounter{proofstep}\ruimte \par 
          \noindent{\bf Step \theproofstep} \hspace{5pt}}
\newcommand{\firststep}{\setcounter{proofstep}{0}\nextstep}

\newcommand{\Ii}{{\bf I}}
\newcommand{\Ni}{{\bf N}}
\newcommand{\Ri}{{\bf R}}
\newcommand{\Ci}{{\bf C}}
\newcommand{\Ti}{{\bf T}}
\newcommand{\Zi}{{\bf Z}}

\newcommand{\proof}{\mbox{\bf Proof} \hspace{5pt}} 
\newcommand{\remark}{\mbox{\bf Remark} \hspace{5pt}}
\newcommand{\ruimte}{\vskip10.0pt plus 4.0pt minus 6.0pt}

\newcommand{\ad}{\mathop{\rm ad}}
\newcommand{\Ad}{\mathop{\rm Ad}}
\newcommand{\RRe}{\mathop{\rm Re}}
\newcommand{\IIm}{\mathop{\rm Im}}
\newcommand{\Tr}{\mathop{\rm Tr}}
\newcommand{\supp}{\mathop{\rm supp}}
\newcommand{\esssup}{\mathop{\rm ess\,sup}}
\newcommand{\Leibniz}{\mathop{\rm Leibniz}}

\hyphenation{groups}
\hyphenation{unitary}

\newcommand{\sun}{\odot}
\newcommand{\hc}{\overline{H}}
\newcommand{\sodp}{ \{ 1,\ldots,d' \} }
\newcommand{\ccig}{C_c^\infty(G)}

\newcommand{\tfrac}[2]{{\textstyle \frac{#1}{#2} }}
\newcommand{\sptilde}{\hspace{0pt} \widetilde{\rule{0pt}{8pt} \hspace{7pt} }}

\newcommand{\cf}{{\cal F}}
\newcommand{\ch}{{\cal H}}
\newcommand{\cl}{{\cal L}}
\newcommand{\cx}{{\cal X}}
\newcommand{\cy}{{\cal Y}}

\newfont{\fontcmrten}{cmr10}
\newcommand{\slbrl}{\mbox{\fontcmrten (}}
\newcommand{\slbrr}{\mbox{\fontcmrten )}}

\newenvironment{remarkn}{\begin{rem} \rm}{\end{rem}}
\renewcommand{\Box}{\raisebox{0.6ex}{\framebox[0.6em]{\rule{0em}{0.6ex}}}}

\newcommand{\idotsint}{\int}

\renewcommand{\thesubsection}{\thesection.\alph{subsection}}

\thispagestyle{empty}

\vspace*{1cm}
\begin{center}
{\Large{\bf Spectral asymptotics of periodic}}  \\[3mm]
{\Large{\bf elliptic operators}}  \\[3mm]

{\large  Ola Bratteli$^1$,}\\
{\large Palle E. T. J{\o}rgensen$^2$}\\  
{\large and}\\
{\large Derek W. Robinson}
 \end{center}

\vspace{2mm}
\begin{center}
\footnotesize Mathematics Research Report No.\ MRR 002--97
\end{center}

\vspace{1.1cm}
\noindent
Centre for Mathematics and its Applications\\
School of Mathematical Sciences\\
Australian National University\\
Canberra, ACT 0200\\
Australia

\vspace{2mm}
\noindent
January 1997

\vspace{2mm}
\vspace*{2mm}

\begin{center}
{\bf Abstract}
\end{center}

\begin{list}{}{\leftmargin=1.8cm \rightmargin=1.8
cm}
\item
We demonstrate that the structure of complex second-order strongly elliptic
operators $H$ on
$\Ri^d$  with coefficients invariant under translation by $\Zi^d$ can be
analyzed through
decomposition in terms of versions $H_z$, $z\in\Ti^d$, of $H$ with
$z$-periodic boundary
conditions acting on $L_2(\Ii^d)$ where $\Ii=[0,1\rangle$.
If the semigroup $S$ generated by $H$ has a H\"older continuous integral
kernel satisfying
Gaussian bounds then the semigroups $S^z$ generated by the $H_z$ have 
kernels with similar
properties and $z\mapsto S^z$ extends to a function on
$\Ci^d\backslash\{0\}$ which is analytic
with respect to  the trace norm.
The sequence of semigroups $S^{(m),z}$ obtained by
rescaling the coefficients of $H_z$ by $c(x)\to c(mx)$ converges in trace
norm to the semigroup
${\widehat S}^z$ generated by the homogenization ${\widehat H}_z$ of $H_z$.
These convergence properties allow asymptotic analysis of the spectrum of $H$.

\end{list}

\vspace{2mm}
\noindent

\begin{enumerate}
\itemindent=-6mm                       
\item Permanent address: 
University of Oslo, P. B. 1053 Blindern, N-0316 Oslo 3, Norway 
\item Permanent address: 
University of Iowa, Iowa City,
IA-52242-1466 USA 
\end{enumerate}

\newpage
\setcounter{page}{1}

\section{Introduction}

We analyze complex, strongly elliptic, periodic operators $H$ on 
 $L_2(\Ri^d)$ in the high frequency limit.
We assume the coefficients of $H$ are invariant under translation 
by the group $\Zi^d$ and
 demonstrate that the semigroup $S$  generated by $H$ decomposes as a
direct integral
of semigroups $S^z$, $z\in\Ti^d$, generated by versions $H_z$ of 
$H$ acting on
$L_2(\Ii^d)$, where $\Ii^d=\Ri^d/\Zi^d$, with $z$-periodic boundary 
conditions.
The decomposition corresponds to a partial Fourier decomposition of 
$L_2(\Ri^d)$ of the
type originally occurring in Bloch wave theory \cite{Blo} which has 
 recently been used 
in wavelet theory (see, for example, \cite{Dau}, pages~109--112). 
It has the advantage that if $S$ has an integral kernel with the
usual continuity and boundedness properties then the $S^z$ have kernels 
which inherit similar
properties. 
Since the $S^z$ act on $L_2(\Ii^d)$ and $\Ii^d$ is bounded  
it  follows that the
$S^z_t$, $t>0$, are trace class operators.
The spectrum of $H$ can then be analyzed by combination of 
decomposition theory and the
spectral theory of the $S^z$ in terms of Bloch bands.

The spectral analysis of the elliptic operator  $H$ can be approached
by homogenization theory (see, for example, \cite{BLP} or \cite{ZKON}).
This corresponds to scaling the period to zero. 
In particular if $H^{(m)}$ denotes the sequence of operators obtained 
from $H$ by rescaling
the coefficients of $H$ by the replacement 
$c_{ij}(x)\to c^{(m)}_{ij}(x)=c_{ij}(mx)$ etc.
then the semigroups $S^{(m)}$ generated by the $H^{(m)}$ 
converge in norm, on each
$L_p(\Ri^d)$-space, to the semigroup $\widehat S$ 
generated by the constant
coefficient homogenization $\widehat H$ of $H$.
This is proved by a variation of methods used in our 
earlier paper \cite{BBJR} with Charles
Batty and the proof is based on  kernel properties.
It has the remarkable consequence that the sequence of 
semigroups $S^{(m),z}$ obtained from
rescaling the $z$-periodic semigroups $S^z$ 
converges in trace norm on $L_2(\Ii^d)$ to the
homogenization of ${\widehat S}^z$ of the $S^z$.
Therefore the spectrum of ${\widehat S}^z$ provides an 
asymptotic approximation to that of
$S^z$.

In the case of pure second-order operators with real coefficients,
$H=-\sum^d_{i,j=1}\partial_ic_{ij}\partial_j$, our results 
are in part motivated by geometry
where the $c_{ij}$ refer to the metric tensor.
If a discrete abelian group $\Gamma$ acts freely on a 
non-compact differentiable manifold
$\cal M$ with compact quotient then Atiyah \cite{Ati} 
and Donnelly \cite{Don} describe the
fibering of elliptic operators $D$ on $\Gamma$-covariant 
bundles over ${\cal M}/\Gamma$.
Trace estimates on the corresponding operators $D^z$ in the 
fibers lead to index
computations in the geometric setting.
Some of our results apply to this setting with only minor modification.

In Section~2 we describe the decomposition theory for 
periodic operators, semigroups and
kernels.
In Section~3 we extend our earlier results on 
homogenization and then in Section~4 we
discuss the general notion of spectral refinement 
in the high frequency limit.
We conclude with some remarks on Schr\"odinger operators.

\section{Periodic decompositions}

In this section we examine  decomposition theory for second-order
 elliptic operators on the complex
Hilbert space $L_2(\Ri^d)$ with periodic coefficients.
Throughout the section we suppose that $H$ is the maximal accretive
 operator associated with the sectorial
form (see, for example, \cite{Kat1} \cite{RS4})
\begin{equation}
h(f)=\sum^{d}_{i,j=1}\,(\partial_if,c_{ij}\partial_jf)+
\sum_{i=1}^d\Big(
({\overline c}_if,\partial_if)+(\partial_if,c'_if)\Big) 
+(f,c_0f)
\label{per2.1}
\end{equation}
where the $\partial_i=\partial/\partial x_i$ denote the 
usual partial derivatives
and the domain of $h$ is $D(h)=L_{2;1}(\Ri^d)=\bigcap_{i=1}^d 
D(\partial_i)$.
The complex-valued coefficients $c_{ij}, c_i,c'_i,c_0\in
 L_\infty(\Ri^d)$ and the matrix of principal
coefficients $C=(c_{ij})$, which is not necessarily
symmetric,  is assumed to satisfy the ellipticity condition
\[
\Re C=(C+C^*)/2\geq \lambda I>0\;\;\;,
\]
in the sense of $d\times d$-matrices over $\Ci^d$, uniformly over $\Ri^d$.
It is this condition which ensures that the form $h$ is sectorial.
The least upper bound $\lambda_C$ of the $\lambda$ satisfying the 
ellipticity condition is called the
ellipticity constant. 

The operator $H$ automatically generates a continuous holomorphic 
semigroup $S$ on $L_2(\Ri^d)$ but there
is little one can deduce about the action of $S$ with no further 
reality, symmetry or smoothness
assumptions on the coefficients. 
It does follow by a perturbation argument \cite{AMT1} \cite{ER20} 
that $S$ leaves $L_2(\Ri^d)\cap
L_p(\Ri^d)$  invariant for $p$ sufficiently close to $2$ and for 
all $p\in[1,\infty]$ if $d=1$ or $d=2$.
Then $S$ extends to a continuous semigroup on the appropriate 
$L_p$-space.
But there are examples \cite{ACT} which show that $S$ does not 
necessarily extend to all the $L_p$-spaces
if $d\geq5$. 
We  first use  periodicity of the coefficients and then kernel bounds to
analyze the action of $S$.

All subsequent estimates depend on $\lambda_C$ and the 
$L_\infty(\Ri^d)$-norms of the coefficients of $H$.
In order to trace the uniformity of the estimates it is 
convenient to introduce ${\cal E}_N$, for each
$N>0$, as the set of elliptic operators $H$ of the 
above type with
\[
\lambda_C^{-1}+\sum^{d}_{i,j=1}\|c_{ij}\|_\infty+
\sum_{i=1}^d\Big(\|c_i\|_\infty+\|c_i'\|_\infty\Big)+
\|c_0\|_\infty\leq N\;\;\;.
\]
In addition we use ${\cal E}^0_N$ to denote the subset 
of ${\cal E}_N$
consisting of the pure second-order operators, i.e., 
those with $c_i=c_i'=c_0=0$.
In fact the magnitude of $\RRe c_0$ is not important 
for uniformity of most of the estimates so we could
use the alternative space  ${\cal E}_N'$ consisting of 
those $H$ for which there is a $\mu_N\in\Ri$
such that $H+\mu_N I\in {\cal E}_N$.

It follows by an elementary estimation that there is a 
$\mu\geq0$ such that
\[
\RRe h(f)\geq-\mu\|f\|_{L_2(\Ri^d)}
\]
for all $f\in D(h)$ where the value of $\mu$ depends only on the 
ellipticity constant and the
$L_\infty(\Ri^d)$-norms of the coefficients. 
Therefore, for each $N>0$ there is a $\mu_N\geq0$ such that 
$\Re H\geq-\mu_N I$ uniformly for all 
$H\in {\cal E}_N$ where $\Re H$ denotes the self adjoint 
operator associated with the closed quadratic form
$\RRe h$. 
Then by the addition of $\mu_N$ to
$c_0$ we may assume $\Re H$ is positive for all $H\in {\cal E}_N$.
The convention $\Re H\geq0$ ensures that the semigroup $S$ 
generated by $H$ is contractive.
The angle of the holomorphy sector of $S$ can then  be 
estimated in terms of the coefficients.
In particular $S$ is holomorphic in the interior of the 
sector $\Delta(\theta)=
\{\zeta:\,|\arg\zeta|\leq\cot^{-1} (\|\Im C\|_\infty/\lambda_C)\,\}$ 
where $\|\Im C\|_\infty$ denotes the
$L_\infty(\Ri^d)$-norm  of the norm of the matrix $\Im C=(C-C^*)/2i$.
Thus for each $N>0$ the semigroups $S$ generated by the $H\in{\cal E}_N$ 
have a common open sector of
holomorphy.

In addition to the general elliptic structure we  assume throughout 
that  the  coefficients of $H$ are
periodic, i.e.,
\begin{equation}
c_{ij}(x+n)=c_{ij}(x),\;\;\;\;c_i(x+n)=c_i(x),\;\;\;\;
{\rm etc.}\label{per2.2}
\end{equation}
for all $x\in \Ri^d$ and $n\in\Zi^d$.
(We have chosen the periods equal to one for simplicity.)
If $U$ denotes the unitary action of $\Ri^d$ by left translations 
on  $L_2(\Ri^d)$, i.e.,
\[
(U(y)f)(x)=f(x-y)
\]
for all $f\in L_2(\Ri^d)$ and $x,y\in\Ri^d$, then $U(x)D(h)=D(h)$ 
for all $x\in\Ri^d$ and the periodicity 
(\ref{per2.2}) of the coefficients  gives the invariance property 
$h(U(n)f)=h(f)$ for all $n\in \Zi^d$.
Hence $U(x)D(H)=D(H)$ and $U(n)H=HU(n)$.
The periodicity of the coefficients of $H$ is reflected by the 
commutation relations
\[
U(n)S_\zeta=S_\zeta U(n)
\]
on $L_2(\Ri^d)$ for all $n\in\Zi^d$ and all $\zeta$ in the 
holomorphy sector of $S$.
Next we examine versions of $H$ and $S$ on $L_2(\Ii^d)$ where
 $\Ii=[0,1\rangle$.

First, introduce  the partial derivatives $\partial_i^z$ on 
$L_2(\Ii^d)$ as the 
skew-adjoint operators of differentiation corresponding to  the 
$z$-periodic boundary conditions, 
$f(u_1,\ldots,1,\ldots, u_d)=z_if(u_1,\ldots,0,\ldots, u_d)$ 
where the $1$ and $0$ are in the $i$-th
position and $z_i\in\Ti$. 
Secondly, define $H_z$  as the maximal accretive operator on 
$L_2(\Ii^d)$ associated with the 
sectorial form
\begin{equation}
h_z(f)=\sum^{d}_{i,j=1}\,(\partial^z_if,c_{ij}\partial^z_jf)+\sum_{i=1}^d
\Big(({\overline c}_if,\partial^z_if)+ (\partial^z_if,c'_if)\Big)
+(f,c_0f)
\label{per2.3}
\end{equation}
where $D(h_z)=\bigcap_{i=1}^d D(\partial^z_i)$.
Repetition of the ellipticity estimates which  gave 
$\Re H\geq-\mu_N I$ then gives $\Re H_z\geq-\mu_N I$
for each $z\in\Ti^d$.
Thus the normalization $\Re H\geq0$ ensures the 
$\Re H_z$ are also positive. 
The $H_z$  are versions of $H$ with $z$-periodic 
boundary conditions and as a consequence have discrete
spectrum. 
\begin{lemma}\label{lper2.0}
The operators $H_z$, $z\in\Ti^d$,
 have compact resolvents.
\end{lemma}
\proof\
First, the real part of $h_z$ satisfies bounds
\[
\RRe h_z(f)\geq\lambda\sum^d_{i=1}
\|\partial_if\|^2_{L_2(\Ii^d)}-\mu \|f\|^2_{L_2(\Ii^d)}
\]
with $\lambda>0$ for all $f\in  D(h_z)$.
Therefore
\[
(1+\mu)I+\Re H_z\geq I+\lambda L_N
\]
where $L_N$ is the version of the 
Laplacian on $L_2(\Ii^d)$ with Neumann boundary conditions.
Since $L_N$ has compact resolvent 
it follows that $\Re H_z$ has compact resolvent.

Secondly, since we have the 
normalization convention $\Re H_z\geq0$ one can represent the resolvent
$(\lambda I+H_z)^{-1}$ as 
\[
(\lambda I+H_z)^{-1}=(\lambda I+\Re H_z)^{-1/2}
(I+C_\lambda)(\lambda I+\Re H_z)^{-1/2}
\]
for all $\lambda>0$ where $C_\lambda$ is bounded 
(see, for example, (3.8) in \cite{Kat2}).
Since $(\lambda I+\Re H_z)^{-1/2}$ is compact it follows that  
$(\lambda I+H_z)^{-1}$ must be compact for
all $\lambda>0$.\hfill\Box
\bigskip

It follows from the positivity and  maximal accretivity that the 
$H_z$ generate a family of
strongly continuous contraction semigroups $S^z$ on $L_2(\Ii^d)$
 and we subsequently demonstrate that the
semigroup $S$ has a decomposition in terms of the $S^z$.  
We show in the next lemma that the family $z\mapsto S^z$ is 
strongly continuous.
Subsequently, in Corollary~\ref{cper2.1}, we establish under 
slightly stronger assumptions that $z\mapsto
S^z$ extends to a function on $\Ci^d\backslash\{0\}$ which is 
analytic with respect to either the
Hilbert--Schmidt norm or the trace norm.
\begin{lemma}\label{lper2.1}
For each $z_0\in \Ti^d$ and $f\in L_2(\Ii^d)$
\[
\lim_{z\to z_0}\|S^z_tf-S^{z_0}_tf\|_{L_2(\Ii^d)}=0
\]
 uniformly for $t$ in any finite  interval of $[0,\infty\rangle$.
\end{lemma}
\proof\
It is convenient for the proof to use the parametrization $z_j=
e^{i\theta_j}$ with $\theta_j\in[-\pi,\pi]$
and to replace the $z$ indices and suffices by $\theta$.

Define the  map $\varphi\in[-\pi,\pi]^d\mapsto V(\varphi)$ into 
unitaries on $L_2(\Ii^d)$ by
\[
(V(\varphi)f)(u)=e^{iu.\varphi}f(u)
\;\;\;.
\]
This map is norm continuous and in particular $\|I-V(\varphi)
\|_{2\to 2}\to0$ as $|\varphi|\to0$.
Moreover, $V(\varphi)D(h_\theta)=D(h_{\theta+\varphi})$.
But
\begin{equation}
\partial_jV(\varphi)=V(\varphi)(\partial_j+i\varphi_j)
\;\;\;.\label{per2.31}
\end{equation}
Therefore
\[
h_{\theta+\varphi}(V(\varphi)f)=h_\theta(f)+p_{\theta,\varphi}(f)=
h_{\theta,\varphi}(f)
\]
for all $f\in D(h_\theta)$ where  $p_{\theta,\varphi}$ is a small 
form perturbation of $h_\theta$
with a relative bound which tends to zero as $|\varphi|\to0$.
This last property is an immediate consequence of (\ref{per2.3}) 
and (\ref{per2.31}).
Thus if $S^{\theta,\varphi}$ denotes the continuous semigroup 
generated by the  maximal accretive
operator $H_{\theta,\varphi}$  associated with the form
$h_{\theta,\varphi}$ one has
\[
V(\varphi)^*S^{\theta+\varphi}_tV(\varphi)=
S_t^{\theta,\varphi}\;\;\;.
\]
and hence  $\|S_t^{\theta,\varphi}\|_{2\to2}\leq 1$ 
for all $\varphi$.
Moreover $ S_t^{\theta,\varphi}$ converges strongly to 
$S^\theta_t$ as $|\varphi|\to 0$ and the
convergence is uniform for $t$ in any finite  interval of 
$[0,\infty\rangle$.
Hence 
\[
(S^{\theta+\varphi}_t-S^{\theta}_t)f=
V(\varphi)(S_t^{\theta,\varphi}-S^\theta_t)V(\varphi)^*f+
(V(\varphi)S_t^{\theta}V(\varphi)^*-S^\theta_t)f\;\;\;.
\]
Now the desired convergence as $|\varphi|\to0$ follows 
by a simple estimate using the norm convergence
of $V(\varphi)\to I$ and the strong convergence of 
$S_t^{\theta,\varphi}\to
S^\theta_t$.\hfill\Box
\bigskip

Next we examine the decomposition of $S$ in terms of the $S^z$.

First define the Zak transform (see \cite{Dau} pages~109--112)  
$Z\colon L_2(\Ri^d)\mapsto 
L_2(\Ti^d\times \Ii^d)$ by
\begin{equation}
(Zf)(z\,,u)=\sum_{n\in\Zi^d}z^nf(u-n)\label{per2.4}
\end{equation}
and  $z^n=z_1^{n_1}z_2^{n_2}\ldots z_d^{n_d}$.
If $f$ has compact support then the sum is finite and 
one calculates straightforwardly that
\begin{eqnarray}
\|Zf\|^2_{L_2(\Ti^d\times\Ii^d)}&=&
(2\pi)^{-d}\int_{\Ti^d}|dz|\int_{\Ii^d}du\,
\Big|\sum_{n\in\Zi^d}z^nf(u-n)\Big|^2
\nonumber\\[5pt]
&=&\sum_{n\in\Zi^d}\int_{\Ii^d}du\,|f(u-n)|^2=
\|f\|^2_{L_2(\Ri^d)}\label{per2.5}
\end{eqnarray}
where $|dz|=|dz_1|\ldots|dz_d|$.
Therefore $Z$ extends by continuity to an isometric map.
But the inverse map  is defined by
\[
(Z^{-1}f)(u-n)=(2\pi)^{-d}\int_{\Ti^d}|dz|
\,{\overline z}^nf(z\,,u)
\]
for $u\in \Ri^d$ and $n\in\Zi^d$
and one again calculates that $Z^{-1}$ is a 
densely defined isometry.
Hence $Z$  extends to a unitary map from 
$L_2(\Ri^d)$ to $L_2(\Ti^d\times \Ii^d)$.
It is often convenient to extend the definition 
(\ref{per2.4}) to all $u\in\Ri^d$ and the resulting
transformation satisfies the periodicity condition
$(Zf)(z\,,x+n)=z^n(Zf)(z\,,x)$ for all $z\in\Ti^d$, 
$x\in\Ri^d$ and $n\in\Zi^d$.

The Zak transform gives a decomposition of $L_2(\Ri^d)$,
\[
L_2(\Ri^d)=(2\pi)^{-d}\int_{\Ti^d}^\oplus|dz|\,L_2(\Ii^d)_z=
L_2(\Ti^d)\otimes L_2(\Ii^d)\;\;\;,
\]
 as a direct integral of copies  of $L_2(\Ii^d)$
indexed by $z\in\Ti^d$ (for a description of the 
formalism of integral decompositions see, for example,
\cite{Dix} or \cite{BR1}). 
We refer to this as the Zak decomposition of $L_2(\Ri^d)$. 
In particular if $f\in L_2(\Ri^d)$  then, by Fubini's theorem, 
for almost all
$z\in\Ti^d$ the function $f_z$  defined by 
\begin{equation}
f_z(u)=(Zf)(z\,,u)\label{per2.6}
\end{equation}
is in $L_2(\Ii^d)$ and (\ref{per2.5}) states that 
\begin{equation}
\|f\|^2_{L_2(\Ri^d)}= (2\pi)^{-d}\int_{\Ti^d}|dz|\,
\|f_z\|^2_{L_2(\Ii^d)}
=\|Zf\|^2_{L_2(\Ti^d)\otimes L_2(\Ii^d)}\label{per2.61}
\;\;\;.
\end{equation}
Note that $ZU(n)Z^*=M_{z^n}$ where $M_{z^n}$ is the 
operator of multiplication by $z^n$ on 
$L_2(\Ti^d)\otimes L_2(\Ii^d)$.
Thus if $A$ is any bounded operator on $L_2(\Ri^d)$ 
such that $AU(n)=U(n)A$ for all $n\in\Zi^d$ then
$M_{z^n}ZAZ^*=ZAZ^*M_{z^n}$ for all $n\in\Zi^d$ and 
hence $ZAZ^*$ commutes with 
$L_\infty(\Ti^d)\otimes I_{L_2(\Ii^d)}$.
Therefore $ZAZ^*$ is a decomposable operator,
\[
ZAZ^*=(2\pi)^{-d}\int_{\Ti^d}^\oplus|dz|\,A(z)\;\;\;,
\]
where the $A(z)$ are bounded operators on $L_2(\Ii^d)$ 
for $z\in\Ti^d$.
If $A$ is an unbounded operator, or form, more care has 
to be taken but we will generally identify 
operators $A$ on $L_2(\Ri^d)$ and $ZAZ^*$ on 
$L_2(\Ti^d)\otimes L_2(\Ii^d)$.
Now we argue that the operator $H$ and the 
semigroup $S$  decompose in this
manner.

Let $f\in C^\infty_c(\Ri^d)\subset D(h)$ and define 
$f_z$ by (\ref{per2.6}).
Since $f$ is differentiable it follows from the 
properties of the Zak transform that $f_z\in D(h_z)$
and $(\partial_if)_z=\partial^z_if_z$  for all $z\in\Ti^d$.
Therefore a straightforward calculation using the 
periodicity of the coefficients gives
\begin{equation}
h(f)= (2\pi)^{-d}\int_{\Ti^d}|dz|\,h_z(f_z)\;\;\;.\label{per2.8}
\end{equation}
But if $D(h)$ is equipped with the norm $f\mapsto 
(\|f\|^2_{L_2(\Ri^d)}+\RRe h(f))^{1/2}$
and the $D(h_z)$ are equipped with analogous norms 
then (\ref{per2.8}) extends by closure to all of $D(h)$,
by use of (\ref{per2.61}).
Thus if $f\in D(h)$ then the family 
$z\in\Ti^d\mapsto f_z$ is $|dz|$-almost  everywhere in $D(h_z)$ and
\[
H=(2\pi)^{-d}\int_{\Ti^d}^{\oplus}|dz|\,H_z
\]
in the sense of direct integral decompositions 
of closed sectorial forms.

One has a similar decomposition of the semigroup.
\begin{thm}\label{tper2.1} Let $S^z$, $z\in\Ti^d$,
 denote the continuous contraction semigroups generated by
the periodic subelliptic operators $H_z$ on $L_2(\Ii^d)$
 with $z$-periodic boundary conditions.
Then the semigroup $S$ generated by $H$ on $L_2(\Ri^d)$ 
has the integral decomposition
\[
S_t=(2\pi)^{-d}\int_{\Ti^d}^{\oplus}|dz|\,S^z_t
\]
corresponding to the Zak decomposition of $L_2(\Ri^d)$.
\end{thm}
\proof\
The proof is relatively straightforward but requires 
several approximation techniques and some standard
measure theoretic arguments.
We sketch the main ideas.

First it is convenient to assume the coefficients of $H$ 
are $C^\infty$-functions. 
This ensures that $C^\infty_c(\Ri^d)$ is a core of $H$.
Then this smoothness assumption is removed by a limiting argument.

Now if $f\in C_c^\infty(\Ri^d)$ one has $(Hf)_z=H_zf_z$ 
and hence
\[
(\lambda I+H_z)^{-1}g_z=((\lambda I+H)^{-1}g)_z
\]
for all $\lambda>0$ and all $g$ in the dense set 
$D_\lambda=(\lambda I+H)(C_c^\infty(\Ri^d))$.
Therefore
\[
\|(\lambda I+H)^{-1}g\|^2_{L_2(\Ri^d)}=
(2\pi)^{-d}\int_{\Ti^d}|dz|\,\|(\lambda I+H_z)^{-1}g_z
\|^2_{L_2(\Ii^d)}
\]
for all $g\in D_\lambda\subset C_c^\infty(\Ri^d)$.
This relation then extends to all $g\in L_2(\Ri^d)$ and 
$(\lambda I+H_z)^{-1}g_z\in L_2(\Ii^d)$ for all
$z$ in a set $\Omega_\lambda$ with $|dz|$-measure equal to one.
Then by iteration and a diagonalization argument
\[
\|( I+t_iH/n)^{-n}g\|^2_{L_2(\Ri^d)}=
(2\pi)^{-d}\int_{\Ti^d}|dz|\,
\|( I+t_iH_z/n)^{-n}g_z\|^2_{L_2(\Ii^d)}
\]
and $( I+t_iH_z/n)^{-n}g_z\in L_2(\Ii^d)$ 
for all rational $t_i$, all positive integers $n$ and all $z$ in
a set $\Omega$ with $|dz|$-measure equal to one.
But the left hand side converges to 
$\|S_{t_i}g\|^2_{L_2(\Ri^d)}$ as $n\to\infty$ and in addition 
$\|( I+t_iH_z/n)^{-n}g_z\|^2_{L_2(\Ii^d)}\to 
\|S^z_{t_i}g_z\|^2_{L_2(\Ii^d)}$ for all $z\in \Omega$.
Therefore, by the Lebesgue dominated convergence theorem,
\[
\|S_{t_i}g\|^2_{L_2(\Ri^d)}=
(2\pi)^{-d}\int_{\Ti^d}|dz|\,
\|S^z_{t_i}g_z\|^2_{L_2(\Ii^d)}
\]
and $S^z_{t_i}g_z\in L_2(\Ii^d)$ for all 
$z\in\Omega$ and all rational $t_i$.
Hence by continuity of $S$ and  the $S^z$ one has
\[
\|S_{t}g\|^2_{L_2(\Ri^d)}=
(2\pi)^{-d}\int_{\Ti^d}|dz|\,\|S^z_{t}g_z\|^2_{L_2(\Ii^d)}
\]
for all $g\in L_2(\Ri^d)$ and 
$S^z_{t}g_z\in L_2(\Ii^d)$ for all $t\geq0$ and all $z$ in
a set $\Omega$ with $|dz|$-measure equal to one.

Finally, if the coefficients of $H$ are only measurable, 
one can approximate $H$ by a sequence of elliptic
operators $H_n$ with $C^\infty$-coefficients obtained by
 regularization of the coefficients of $H$.
Specifically the coefficients of $H$ are replaced by
\[
c^{(n)}_{ij}(x)=n^d\int_{\Ri^d}dy\,\tau(ny)c_{ij}(x-y)
\]
etc.\ where $\tau$ is a positive $C^\infty$-function with integral one.
Then the ellipticity constants of the $H_n$ are bounded below 
by the ellipticity constant of $H$ and the
$L_\infty$-norms of the regularized coefficients are bounded 
above by the $L_\infty$-norms of the
unregularized coefficients.
Thus if $H\in{\cal E}_N$ then $H_n\in{\cal E}_N$.
Moreover, the sequence of semigroups $S^{(n)}$ generated by 
the $H_n$ converges strongly to $S$ on
$L_2(\Ri^d)$ (see \cite{ER15}, Proposition 2.6).
In particular
\[
\lim_{n\to\infty}\|S^{(n)}_tf-S_tf\|_{L_2(\Ri^d)}=0
\]
for all $f\in L_2(\Ri^d)$ uniformly for $t$ in finite 
intervals of $[0,\infty\rangle$.
But the same regularization procedure can be applied to 
the $H_z$ and one concludes that one also has
strong convergence of the corresponding semigroups 
$S^{(n),z}$ to the $S^z$.
All the estimates required in the regularization argument 
are independent of the particular choice of
boundary condition.
Therefore arguing as before one can deduce that
\[
\|S^{(n)}_{t}g\|^2_{L_2(\Ri^d)}=
(2\pi)^{-d}\int_{\Ti^d}|dz|\,\|S^{(n),z}_{t}g_z\|^2_{L_2(\Ii^d)}
\]
for all $g\in L_2(\Ri^d)$, $n\in\Ni$ and $t>0$. 
Moreover, $S^{(n),z}_{t}g_z\in L_2(\Ii^d)$ for all $t\geq0$, 
all $n\in\Ni$ and all $z$ in
a set $\Omega$ with $|dz|$-measure equal to one.
Taking the limit $n\to\infty$ then gives the required 
decomposition of $S$.\hfill\Box
\bigskip

Next we assume that the action of the semigroup $S$ is given by a 
H\"older continuous integral kernel $K$
satisfying Gaussian upper bounds uniformly for $H$ in each ${\cal E}_N$.
Specifically we assume that for each $N>0$ the semigroup
$S$ generated by  $H\in{\cal E}_N$ has a
kernel $K$ such that
\[
(S_\zeta f)(x)=\int_{\Ri^d}dy\,K_\zeta (x\,;y)f(y)\;\;\;,
\]
for all $f\in L_2(\Ri^d)$ and all $\zeta$ in the sector of holomorphy of 
$S$ satisfying the following
properties:
\begin{enumerate}
\item
there exist $a,b>0$, $\omega\geq0$ and 
$\theta\in\langle0,\pi/2\rangle$ such that
\begin{equation}
|K_\zeta(x\,;y)|\leq a\,|\zeta|^{-d/2}e^{\omega |\zeta|}
e^{-b|x-y|^2/|\zeta|}\label{per2.13}
\end{equation}
for all $x,y\in\Ri^d$ and $\zeta\in\Delta(\theta)$,
\item
in addition there exists $\nu\in\langle0,1]$ such that 
\begin{equation}
|K_\zeta(x-x'\,;y-y')-K_\zeta(x\,;y)|\leq a\,
|\zeta|^{-d/2}e^{\omega
|\zeta|}\bigg({{|x'|+|y'|}\over{|\zeta|^{1/2}}}\bigg)^\nu 
e^{-b|x-y|^2/|\zeta|}\label{per2.14}
\end{equation}
for all $x,x',y,y'\in\Ri^d$ and all $\zeta\in \Delta(\theta)$  
with $|x'|+|y'|\leq
|\zeta|^{1/2}$,
\end{enumerate}
with both bounds uniform for $H\in{\cal E}_N$.
The uniformity of these bounds will be often used in the sequel.

The $K_\zeta$  satisfy  semigroup composition properties which reflect the
corresponding properties of $S$.
Moreover the periodicity of the coefficients is reflected in 
the periodicity of the kernel,
\[
K_\zeta(x\,;y+n)=K_\zeta(x-n\,;y)
\]
for all $x,y\in\Ri^d$ and $n\in\Zi^d$.
The adjoint semigroup $S^*$ has a similar kernel $K^*$ related to $K$ by
$K^*_t(x\,;y)={\overline{K_t(y\,;x)}}$.

One immediate consequence of the Gaussian bounds on the  kernel is that $S$
leaves $L_2(\Ri^d)\cap L_p(\Ri^d)$ invariant for each $p\in[1,\infty]$ and
hence extends to a continuous semigroup, also denoted by $S$, 
on each of the spaces
$L_p(\Ri^d)$.
A less obvious implication is the  continuity estimate
\begin{equation}
|K_{t-t'}(x\,;y)-K_t(x\,;y)|\leq a'\,|t'|\,t^{-1-d/2}
e^{\omega t} e^{-b'|x-y|^2/|t|}\label{per2.15}
\end{equation}
which holds for all $x,y\in\Ri^d$ and $t>0$ with $|t'|\leq t/2$.
These bounds follow from the Duhamel formula
\[
K_t(x\,;y)-K_{t-t'}(x\,;y)=\int^t_{t-t'}ds\,{{d}\over{ds}}K_s(x\,;y)
\]
and the Cauchy representation
\[
K_s(x\,;y)=(2\pi i)^{-1}\int_{C_r(s)}d\zeta\,
{{K_\zeta(x\,;y)}\over{s-\zeta}}
\]
where the integral is over a circle $C_r(s)$ of radius $r$ centred at $s$.
The estimates are automatically  uniform for $H\in{\cal E}_N$.

The existence of a H\"older continuous Gaussian kernel, 
in the above sense, is not an automatic consequence
of our assumptions but it does follow from some weak additional hypotheses, 
e.g., if one of the following
three conditions is satisfied,
\begin{enumerate}
\item
if $d=1$ or $2$,
\item
if the principal coefficients $c_{ij}$ are real,
\item
if $d\geq3$ and the  $c_{ij}$ are complex but uniformly continuous.
\end{enumerate}
On the other hand no such kernel exists in general if $d\geq 5$.

Results of this nature for low dimensions are well known but 
recent discussions of Gaussian bounds of the
above type have been given in \cite{AMT1} and \cite{ER20}.
The results for real coefficients are a variation of classic 
work of Nash \cite{Nash} and De Giorgi
\cite{DG}. 
A more recent description of the Nash approach is given in \cite{FaS} 
and De Giorgi's
approach is described in \cite{Gia}.
The situation for complex, uniformly continuous, $c_{ij}$ is 
covered in \cite{Aus1} and \cite{ER15}.
The pathologies of high dimensions are related to those of 
systems of elliptic operators.
A variety of counterexamples in the latter setting is given in \cite{Gia1}.
The failure of Gaussian bounds for $d\geq 5$ is discussed in \cite{ACT}.
Note that if one has the bounds (\ref{per2.13}) and (\ref{per2.14}) 
for complex coefficients and real
$\zeta$ then they follow for complex $\zeta$ in a small sector by `rotation': 
the replacement $t\to
te^{i\theta}$ corresponds to the replacement $c_{ij}\to c_{ij}e^{i\theta}$ etc.

The kernel properties transfer from $S$ to the $S^z$.

\begin{thm}\label{tper2.2} Assume the semigroup $S$ has a 
H\"older continuous kernel $K$ satisfying Gaussian
bounds.
Then each $S^z$ has an integral kernel $K^z$ given by the Zak transform,
\begin{equation}
K^z_\zeta(u\,;v)=\sum_{n\in\Zi^d}z^nK_\zeta(u\,;v+n)=
\sum_{n\in\Zi^d}z^nK_\zeta(u-n\,;v)\label{per2.19}
\end{equation}
of  $K$.

The $K^z$ are jointly  H\"older continuous and satisfy bounds
\begin{eqnarray*}
|K^z_\zeta(u\,;v)|&\leq&  a\,(1\wedge |\zeta|)^{-d/2}
e^{\omega |\zeta|}e^{-b|u-v|^2/|\zeta|}\\[5pt]
|K_\zeta(u-u'\,;v-v')-K_\zeta(u\,;v)|&\leq& a\,
(1\wedge |\zeta|)^{-d/2}e^{\omega
|\zeta|}\bigg({{|u'|+|v'|}\over{|\zeta|^{1/2}}}\bigg)^\nu e^{-b|u-v|^2/|\zeta|}
\end{eqnarray*}
for all $u,v,u',v'\in \Ii^d$  and $\zeta\in \Delta(\theta)$ with 
$|u'|+|v'|\leq |\zeta|^{1/2}$ and for all
$z\in\Ti^d$, uniformly for $H\in{\cal E}_N$.
\end{thm}
\proof\
We give the proof for $\zeta=t>0$.

First define functions $K^z$ by (\ref{per2.19}).
Then the $K^z$ are jointly H\"older continuous and satisfy 
the Gaussian bounds on $\Ii^d$ as a simple
consequence of the definition and the properties of $K$.
Therefore one can define bounded operators ${\widetilde S}_t^z$ on  
$L_2(\Ii^d)$, or on $L_p(\Ii^d)$, by
\[
({\widetilde S}^z_tf)(u)=\int_{\Ii^d}dv\,K^z_t(u\,;v)f(v)\;\;\;.
\] 
But
\begin{eqnarray*}
\int_{\Ii^d}dv\,K^z_s(u\,;v)K^z_t(v\,;w)&=&\sum_{n,m\in\Zi^d}
z^{n+m}\int_{\Ii^d}dv\,K_s(u\,;v+n)K_t(v+n\,;w+n+m)\\[5pt]
&=&\sum_{p\in\Zi^d}
z^p\int_{\Ri^d}dx\,K^z_s(u\,;x)K_t(x\,;w+p)=K^z_{s+t}(u\,;w)
\end{eqnarray*}
where the last identification uses the semigroup property of the kernel $K$.
Hence the
${\widetilde S}^z_t$ form a semigroup
${\widetilde S}^z$ with continuity properties similar to those of $S$. 
Thus it suffices to prove that ${\widetilde S}^z_t=S^z_t$

Let $f\in C_c(\Ri^d)$ then $f_z\in L_2(\Ii^d)$ and
\begin{eqnarray*}
({\widetilde S}^z_tf_z)(u)&=&\sum_{n\in\Zi^d}z^n\int_{\Ii^d}dv\,
K^z_t(u\,;v)f(v-n)\\[5pt]
&=&\sum_{n,m\in\Zi^d}z^{n+m}\int_{\Ii^d}dv\,K_t(u\,;v+m)f(v-n)\\[5pt]
&=&\sum_{n,p\in\Zi^d}z^p\int_{\Ii^d}dv\,K_t(u-p\,;v-n)f(v-n)=
(S_tf)_z(u)\;\;\;.
\end{eqnarray*}
Therefore, from (\ref{per2.61}), one obtains
\[
\|S_tf\|^2_{L_2(\Ri^d)}=(2\pi)^{-d}\int_{\Ti^d}|dz|\,
\|(S_tf)_z\|^2_{L_2(\Ii^d)}
=(2\pi)^{-d}\int_{\Ti^d}|dz|\,\|{\widetilde S}^z_tf_z\|^2_{L_2(\Ii^d)}
\]
and this relation extends to all $f\in L_2(\Ri^d)$ by closure.
But this means that $S$ has a direct integral decomposition in 
terms of the semigroups
${\widetilde S}^z$. Then by Theorem \ref{tper2.1} 
it follows that ${\widetilde S}_t^z=S^z_t$ for
$|dz|$-almost all $z$.
Finally the equality of the semigroups follows 
from the continuity of $z\mapsto S^z_t$ given by Lemma
\ref{lper2.1} and the continuity of 
$z\mapsto {\widetilde S}^z_t$ which is an easy consequence of the
definition of ${\widetilde S}^z$ and the bounds on the kernel $K$. \hfill\Box
\bigskip

The kernel bounds imply that the operators 
$S^z_\zeta$, which are compact as a consequence of
Lemma~\ref{lper2.0}, are in fact Hilbert--Schmidt.
Then the semigroup property implies they are trace class.
Let $\|\cdot\|_{Tr}$ and $\|\cdot\|_{HS}$ denote 
the trace norm and  Hilbert--Schmidt norm respectively.
\begin{lemma}\label{lper2.2}
Let $S=\{S_t\}_{t\geq 0}$  be a semigroup of  
bounded operators on a Hilbert space.
Then
\[
\|S_t\|_{Tr}\leq \|S_{t/2}\|_{HS}^2
\]
for all $t>0$.
Moreover if $S$ and $T$ are two such semigroups then
\[
\|S_t-T_t\|_{Tr}\leq(\|S_{t/2}\|_{HS}+
\|T_{t/2}\|_{HS})\|S_{t/2}-T_{t/2}\|_{HS}
\]
for all $t>0$.
\end{lemma}
\proof\
The first statement of the lemma means that if $S_t$ is 
Hilbert--Schmidt for all $t>0$
then it is also trace class and the norm bounds are valid.
The second statement is interpreted in a similar manner.

The bounds follow from the identities
$S_t=S_{t/2}S_{t/2}$ and
\[
S_t-T_t=S_{t/2}(S_{t/2}-T_{t/2})+(S_{t/2}-T_{t/2})T_{t/2}
\]
together with the observation that
\[
\|AB\|_{Tr}\leq \|A\|_{HS}\|B\|_{HS}
\]
for any pair of Hilbert--Schmidt operators.
\hfill\Box
\bigskip

Now one can estimate the trace norm as follows.
\begin{eqnarray}
\|S^z_{2\zeta}\|_{Tr}&\leq&\|S^z_\zeta\|^2_{HS}=
\int_{\Ii^d}du\int_{\Ii^d}dv\,|K^z_\zeta(u\,;v)|^2\nonumber\\[5pt]
&\leq&\sum_{n\in\Zi^d}\int_{\Ii^d}du\sup_{v\in\Ii^d}|K_\zeta(u-n\,;v)|
\sum_{m\in\Zi^d}\int_{\Ii^d}dv
\sup_{u\in\Ii^d}|K_\zeta(u\,;v+m)|\nonumber\\[5pt]
&=&\int_{\Ri^d}dx\sup_{v\in\Ii^d}|K_\zeta(x\,;v)|
\int_{\Ri^d}dy\sup_{u\in\Ii^d}|K_\zeta(u\,;y)|
\leq a'\,(1\wedge|\zeta|)^{-d}e^{2\omega|\zeta|}\label{eper2.20}
\end{eqnarray}
where the second step uses the periodicity of the kernel, 
the third step follows by estimation with the
Gaussian bounds and
the singularity in $\zeta$ comes from the diagonal 
contributions $x=v$ and $u=y$ with $\omega$  the same
parameter as occurs in (\ref{per2.13}).

This method of estimation  gives the following substantial improvement of
Lemma~\ref{lper2.1}.
\begin{cor}\label{cper2.1}
The family $z\mapsto S^z$ extends to a function on 
$\Ci^d\backslash\{0\}$ which is analytic in the
Hilbert--Schmidt norm, or the trace norm, on $L_2(\Ii^d)$.
\end{cor}
\proof\
The extension of $S^z$ is defined by extension of the
 power series expansion (\ref{per2.19}) of the kernel
$K^z$.
The foregoing estimates then adapt to show that the 
extended $S^z$ are Hilbert--Schmidt operators.
But then for $1/R\leq|z|,|z_0|\leq R$ one has
\begin{eqnarray*}
\|S^z_\zeta-S^{z_0}_\zeta\|^2_{HS}&\leq&
\int_{\Ii^d}du\int_{\Ii^d}dv\,
\bigg(\sum_{n\in\Zi^d}|z^n-z_0^n|\,|K_\zeta(u-n\,;v)|\bigg)^2\\[5pt]
&\leq&|z-z_0|^2\sum_{n\in\Zi^d}nR^n\int_{\Ii^d}du
\sup_{v\in\Ii^d}|K_\zeta(u-n\,;v)|\\
&&\hspace{5cm}{}\cdot\sum_{m\in\Zi^d}mR^m
\int_{\Ii^d}dv\sup_{u\in\Ii^d}|K_\zeta(u\,;v+m)|\\[5pt]
&\leq& a'\,|z-z_0|^2(1\wedge|\zeta|)^{-d}e^{\omega'|\zeta|}
\end{eqnarray*}
where the last estimate relies on the Gaussian bounds.
This establishes the analyticity with respect 
to the Hilbert--Schmidt norm.
The analyticity with respect to the trace norm $\|\cdot\|_{Tr}$ 
is a consequence of the preceding
lemma.
The proof of the corollary is completed by setting 
$S=S^z$ and $T=S^{z_0}$ in the lemma and using the
Hilbert--Schmidt estimates.
\hfill\Box

\bigskip

Finally note the kernel $K$ is pointwise positive 
if and only if the coefficients of $H$ are real.
But even in this situation the $K^z$ are not positive 
except in the purely periodic case $z=1$.
In fact the kernels are complex for $z\neq\pm1$.

\section{Asymptotic properties}

In this section we examine asymptotic properties of the periodic 
second-order operators $H$ considered
in Section~2. 
First we reformulate some of the results obtained in \cite{BBJR} for
periodic operators on $\Ri^d$ and then we use the decomposition theory to
deduce asymptotic properties of the
operators on $\Ii^d$. 
Since we did not consider complex operators nor operators with lower order
terms the
proofs of \cite{BBJR} require some adaptation. 
Throughout this section we assume the semigroups $S$ generated by the elliptic
operators $H$ have H\"older continuous kernels satisfying the Gaussian bounds 
(\ref{per2.13}) and
(\ref{per2.14}) uniformly for $H\in{\cal E}_N$ and for each $N>0$.

If $m$ is a positive integer we define the rescaling $H^{(m)}$ of $H$ 
by the replacement $c_{ij}(x)\to
c^{(m)}_{ij}(x)=c_{ij}(mx)$ and $c_i(x)\to
c^{(m)}_i(x)=c_i(mx)$ etc.
Then each $H^{(m)}$ is still periodic and in fact has periods $1/m$.
Next define the homogenization $\widehat H$ of $H$ as the elliptic operator
with constant 
coefficients ${\widehat c}_{ij}$, ${\widehat c}_{i}$, etc. defined as  follows.
First
\[
{\widehat c}_{ij}=\int_{\Ii^d}du\,c_{ij}(u)-
\sum^d_{k,l=1}(c_{ik},X_{kl}c_{lj})
\]
where $X_{kl}$ are the bounded operators associated with the forms
\[
x_{kl}(f)=(\partial_kf,H_1^{-1}\partial_lf)
\] 
on $L_2(\Ii^d)$ with $\partial_k=\partial^z_k|_{z=1}$ and
$H_1=H_z|_{z=1}$, i.e., with $z=(1,\ldots,1)$.
(The operators $X_{kl}$ are bounded because, by spectral theory,
$H_1^{-1}$ is bounded on the orthogonal complement of the 
constant functions in $L_2(\Ii^d)$.)
Secondly,
\[
{\widehat c}_{i}=\int_{\Ii^d}du\,(c'_i(u)+c_{i}(u))-
\sum^d_{k,l=1}(c_{k},X_{kl}c_{li})\;\;\;.
\]
Thirdly,
\[
{\widehat c}_0=\int_{\Ii^d}dx\,c_0(x)\;\;\;.
\]
This definition of $\widehat H$ coincides with that of \cite{BLP}, pages~16
and 184,
 and in the case of pure
second-order operators with that of \cite{BBJR}. 
Note that if $C=(c_{ij})$ is hermitian then ${\widehat
C}=({\widehat c}_{ij})$ is also hermitian but symmetry of $C$ does not
necessarily 
imply symmetry of
${\widehat C}$. 
Nevertheless the homogenized principal coefficients may be symmetrized
because the
anti-symmetric part gives no contribution to
$\widehat H$.

The $\widehat H$ occur as  limits  of the rescaled versions $H^{(m)}$ of $H$.
The basic argument of homogenization theory \cite{BLP} establishes  
local weak convergence of the
weak solutions $u^{(m)}$ of $(\lambda I+H^{(m)})u^{(m)}=f$ to the 
solution ${\hat u}$ of 
$(\lambda I+{\widehat H}){\hat u}=f$.
This argument is based on the following proposition.
\begin{prop}\label{pper3.0}
There is a $\lambda_0>0$ such that for each $f\in L_2(\Ri^d)$ and 
$\lambda\geq\lambda_0$ the equation
\[
(\lambda I+H^{(m)})u^{(m)}=f
\]
has a unique solution $u^{(m)}\in L_2(\Ri^d)$ satisfying bounds
\[
\|u^{(m)}\|_{L_2(\Ri^d)}+\sum^d_{i=1}
\|\partial_iu^{(m)}\|_{L_2(\Ri^d)}\leq c
\]
uniformly in $m$.
\end{prop}
\proof\
The proof of the proposition in the special case $c'_i=c_i=0$ is given in 
\cite{BLP} but the general
case then follows by the same arguments once one has the {\it a priori} estimate
\[
\|(\lambda I+H^{(m)})f\|^2_{L_2(\Ri^d)}\geq
c\,\Big(\|f\|^2_{L_2(\Ri^d)}+\sum^d_{i=1}
\|\partial_if\|^2_{L_2(\Ri^d)}\Big)
\]
 valid for all sufficiently large $\lambda$ with $c>0$ independent of $m$.
But this is the $L_2$-G\aa rding inequality which follows from the
observation that
\[
\|(\lambda I+H^{(m)})f\|^2_{L_2(\Ri^d)}\geq \lambda^2
\|f\|^2_{L_2(\Ri^d)}+2\lambda\RRe(f,H^{(m)}f)
\]
by the usual ellipticity estimates.
Since these estimates depend only on the ellipticity constant 
and the $L_\infty$-norms of the coefficients
they are uniform in $m$.\hfill\Box
\bigskip

Using the proposition and the arguments of \cite{BLP}, 
Sections~I.2 and I.13, one then proves that the
$u^{(m)}$ of the proposition converge locally to the  solution ${\hat u}$ of 
$(\lambda I+{\widehat H}){\hat u}=f$. 
Then one can use general arguments of functional analysis to convert 
this into strong resolvent
convergence on $L_2(\Ri^d)$ (\cite{BBJR}, page 137).
This in turn implies that the semigroups $S^{(m)}$ generated by the 
$H^{(m)}$ are strongly convergent on
$L_2(\Ri^d)$ to the semigroup $\widehat S$ generated by  $\widehat H$.
In fact we demonstrate below that the $H^{(m)}$ are norm resolvent convergent to
$\widehat H$ and hence the $S^{(m)}$ are norm convergent to  
$\widehat S$ on each of the
$L_p(\Ri^d)$-spaces (see Corollary~\ref{cper3.1}).

The basic convergence property of the kernels of the rescaled operators 
is a version of Theorem~III.4 of
\cite{ZKON}.
\begin{lemma}\label{lper3.2} 
If $K^{(m)}$ is the kernel corresponding to $H^{(m)}$ and $\widehat K$  
the kernel corresponding to
$\widehat H$ then
\[
\lim_{m\to\infty}\sup_{|x|^2+|y|^2\leq vt}
|K^{(m)}_t(x\,;y) -{\widehat K}_t(x\,;y)|=0
\]
for all $v>0$ and $t>0$.
\end{lemma}

The lemma is established in the course of the proof of Theorem~III.4 of
\cite{ZKON}.
The latter proof  has to be modified to
cover the current context of complex coefficients but there is only 
one small change needed.
The first part of the proof uses the semigroup convergence discussed 
prior to the lemma.
The second part of the proof, page 137, invokes the result of Nash 
which gives equicontinuity of the
kernel, which is then used in combination with the Arzela--Ascoli 
theorem to deduce a compactness property. 
But the required equicontinuity now follows from the  bounds 
(\ref{per2.14}) together with the 
bounds (\ref{per2.15}).
\bigskip

Lemma \ref{lper3.2} together with the periodicity and the 
Gaussian bounds now give a semigroup convergence
result which was contained in \cite{BBJR} for the case of 
real symmetric coefficients although it was not
explicitly stated.

\begin{prop}\label{pper3.1}
Let $S^{(m)}$ and $\widehat S$ denote the semigroups generated 
by the rescaled versions $H^{(m)}$ of $H$
and the homogenized operator $\widehat H$.
Then
\[
\lim_{m\to\infty}\|S^{(m)}_t-{\widehat S}_t\|_{L_p(\Ri^d)\to L_p(\Ri^d)}
=0
\]
for all $p\in[1,\infty]$, uniformly for $t$ in  
compact intervals of $\langle 0,\infty\rangle$.
\end{prop}
\proof\
It suffices to prove the convergence on $L_\infty(\Ri^d)$ for the semigroups and
their adjoints since the general result then follows by interpolation.
But as the adjoint semigroups $S^*$  are of the same type as $S$ 
it suffices to prove the $L_\infty(\Ri^d)$
convergence. 
This, however, is equivalent to proving that
\begin{equation}
\lim_{m\to\infty}\sup_{x\in\Ri^d}\int_{\Ri^d}dy\,
|K^{(m)}_t(x\,;y) -{\widehat K}_t(x\,;y)|=0
\label{eper3.3}
\end{equation}
with the correct uniformity in $t$.
But
\[
\int_{\Ri^d}dy\,|K^{(m)}_t(x\,;y) -{\widehat K}_t(x\,;y)|=
\int_{\Ri^d}dy\,|K^{(m)}_t(x-n\,;y) -{\widehat K}_t(x-n\,;y)|
\]
for each $n\in\Zi^d$ by periodicity of the kernels.
Therefore the supremum in (\ref{eper3.3}) can be restricted to 
$x$ with $|x_i|\leq 1$.
Moreover, the Gaussian bounds on $K^{(m)}$, which are uniform in 
$m$ by assumption, imply that
\[
\sup_{|x_i|\leq1}\int_{|y|^2\geq vt}dy\,
|K^{(m)}_t(x\,;y)|\leq a\,e^{-bv^2} 
\]
with similar bounds on the integral with $\widehat K$.
Thus for $\varepsilon, t>0$ one may choose $v$ sufficiently large that 
\[
\int_{\Ri^d}dy\,|K^{(m)}_t(x\,;y) -{\widehat K}_t(x\,;y)|\leq
\varepsilon+\int_{|x|^2+|y|^2\leq vt}dy\,
|K^{(m)}_t(x\,;y) -{\widehat K}_t(x\,;y)|
\;\;\;.
\]
Then the statement of the proposition follows from 
Lemma~\ref{lper3.2}.\hfill\Box
\bigskip

One immediately deduces that the resolvents of the 
$H^{(m)}$ converge in norm.

\begin{cor}\label{cper3.1}
The sequence $H^{(m)}$ converges to $\widehat H$ in the 
norm resolvent sense on $L_p(\Ri^d)$ for each
$p\in[1,\infty]$.
\end{cor}

This follows from  Proposition~\ref{pper3.1} by Laplace transformation.

\bigskip

Next we use these estimates to examine convergence of 
the semigroups $S^{(m),z}$ corresponding to the
rescaling $H^{(m)}_z$ of the operator $H_z$ with 
$z$-periodic boundary conditions.
\begin{thm}\label{tper3.1}
Let ${\widehat S}^z$ denote the semigroup generated by 
the homogenized operator $\widehat H$ on
$L_2(\Ii^d)$ with $z$-periodic boundary conditions.
Then
\[
\lim_{m\to\infty}\sup_{z\in\Ti^d}
\|S^{(m),z}_t-{\widehat S}^z_t\|_{HS}=0
\]
and 
\[
\lim_{m\to\infty}\sup_{z\in\Ti^d}
\|S^{(m),z}_t-{\widehat S}^z_t\|_{Tr}=0
\]
uniformly for $t$ in compact intervals of 
$\langle0,\infty\rangle$.
\end{thm}
\proof\
Introduce $D^{(m)}_t$ by
\[
D^{(m)}_t(x\,;y)=|K^{(m)}_t(x\,;y) -{\widehat K}_t(x\,;y)|
\]
for $x,y\in\Ri^d$.
Then
\begin{eqnarray*}
\|S^{(m),z}_t-{\widehat S}^z_t\|^2_{HS}&\leq&
\int_{\Ii^d}du\int_{\Ii^d}dv\,
|\sum_{n\in\Zi^d}D^{(m)}_t(u-n\,;v)|^2\\[5pt]
&\leq&\int_{\Ri^d}dx\sup_{v\in\Ii^d}D^{(m)}_t(x\,;v)
\int_{\Ri^d}dy\sup_{u\in\Ii^d}D^{(m)}_t(u\,;y)
\end{eqnarray*}
by the estimation procedure used to deduce (\ref{eper2.20}).
But the $D^{(m)}_t$ satisfy Gaussian bounds uniformly in $m$.
Hence the contribution to the integrals for $|x|^2\geq vt$, 
or $|y|^2\geq vt$, can be made arbitrarily
small, uniformly in $m$, by choosing $v$ sufficiently large.
The remaining contributions tend to zero, however, as $m\to\infty$ by 
Lemma~\ref{lper3.2}.
Therefore one has the Hilbert--Schmidt convergence in the uniform sense.

Finally it follows from Lemma~\ref{lper2.2} that
\[
\|S^{(m),z}_t-{\widehat S}^z_t\|_{Tr}\leq
(\|S^{(m),z}_{t/2}\|_{HS}+ \|{\widehat S}^z_{t/2}\|_{HS})
\|S^{(m),z}_{t/2}-{\widehat S}^z_{t/2}\|_{HS}
\]
but the Hilbert--Schmidt norms of the $S^{(m),z}_{t/2}$ and 
${\widehat S}^z_{t/2}$ are bounded uniformly in
$m$ and $z$ by the kernel estimate that gave (\ref{eper2.20}).
\hfill\Box

\bigskip

The trace norm convergence of the semigroups $S^{(m),z}$ to 
${\widehat S}^z$ is equivalent to trace norm
convergence of the resolvents $(\lambda I+H^{(m),z})^{-1}$ to 
$(\lambda I+{\widehat H}^z)^{-1}$ uniformly
for $\lambda$ in compact intervals of $\langle0,\infty\rangle$.
If we specialize to self-adjoint operators, i.e., if we assume
$C=(c_{ij})$ is hermitian,
$c_i'=-{\overline c}_i$ and $c_0$ is real, then one can 
immediately deduce strong statements on convergence
of the spectrum of the operators.
Let $\lambda_n(m,z)$ denote the eigenvalues of $H^{(m),z}$ 
in increasing order repeated according
to multiplicity and ${\widehat \lambda}_n(z)$ the 
corresponding list for ${\widehat H}^z$ (see
(\ref{eper4.37}), below). 
The trace norm convergence implies that the corresponding finite-dimensional
eigenprojections are norm convergent to the eigenprojections of the 
homogenized operator and the eigenvalues
are pointwise convergent, i.e., 
$\lambda_n(m,z)\to {\widehat \lambda}_n(z)$ as $m\to\infty$.
But the trace norm estimate gives a uniform bound on 
the eigenvalue convergence.
One has
\begin{equation}
\sum_{n=0}\Big|e^{-t\lambda_n(m,z)}-
e^{-t{\widehat \lambda}_n(z)}\Big|\leq 
\|S^{(m),z}_t-{\widehat S}^z_t\|_{Tr}\label{eper3.4}
\end{equation}
for all $t>0$ (see \cite{Pow}, Section~5).

If $c'_i=c_i=c_0=0$ the foregoing results can be
 rephrased with the aid of scaling.
The advantage of the pure second-order operators is that 
they are homogeneous of order two under dilations.
Therefore  by change of variable one deduces the scaling property
\begin{equation}
K_t(x\,;y)=m^{-d}K_{m^{-2}t}^{(m)}(m^{-1}x\,;m^{-1}y)\label{eper3.1}
\end{equation} 
of the kernels associated with $H$ and the $H^{(m)}$.
One immediate implication of this scaling is that (\ref{per2.13}) 
is valid with $\omega=0$.
This follows because (\ref{per2.13}) is valid for $K^{(m)}$ 
uniformly in $m$ and hence
\[
|K_t(x\,;y)|\leq  a\,t^{-d/2}e^{\omega m^{-2}t}e^{-b|x-y|^2/t}
\]
for all $m$.
In the limit $m\to\infty$ one obtains the bounds with  $\omega=0$.
Similarly (\ref{per2.14}) is established with $\omega=0$ and as a 
consequence one can take $\omega=0$ in
(\ref{per2.15}) and in the Gaussian estimates on $K^z$.
Next note that since the lower order terms are zero the kernel 
$\widehat K$ obeys the scaling relation
\begin{equation}
{\widehat K}_t(x\,;y)=m^{-d}{\widehat K}_{m^{-2}t}
(m^{-1}x\,;m^{-1}y)\;\;\;.\label{eper3.2}
\end{equation}
Therefore the statement of the Lemma~\ref{lper3.2}  is equivalent to
\[
\lim_{t\to\infty}\sup_{|x|^2+|y|^2\leq vt}t^{d/2}
|K_t(x\,;y) -{\widehat K}_t(x\,;y)|=0
\]
for all $v>0$.
But this is the conclusion of Theorem~III.4 of \cite{ZKON} 
(which was rephrased as Proposition~4.2 of
\cite{BBJR} with the unfortunate omission of the factor $t^{d/2}$).
Finally observe that dilations are isometrically 
implemented on $L_p(\Ri^d)$ by the operators
\[
(V(m)f)(x)=m^{-d/p}f(m^{-1}x)
\]
for all $f\in L_p(\Ri^d)$.
Then one calculates from the homogeneity of $H$ under scaling 
and the definition of $H^{(m)}$ that
\[
V(m)^*(S_t-{\widehat S}_t)V(m)=
S^{(m)}_{m^{-2}t}-{\widehat S}_{m^{-2}t}\;\;\;.
\]
Consequently
\[
\|S^{(m)}_t-{\widehat S}_t\|_{L_p(\Ri^d)\to L_p(\Ri^d)}
=\|S_{m^2t}-{\widehat S}_{m^2t}\|_{L_p(\Ri^d)\to L_p(\Ri^d)}
\;\;\;.
\]
Thus the statement of Proposition~\ref{pper3.1} is equivalent 
with the following asymptotic identification.
\begin{cor}\label{cper3.2} If $c'_i=c_i=c_0=0$  then
\[
\lim_{t\to\infty}\|S_t-{\widehat S}_t\|_{L_p(\Ri^d)\to L_p(\Ri^d)}=0
\]
for all $p\in[1,\infty]$.
\end{cor}
This was the principal conclusion of \cite{BBJR}.
It is equivalent to the convergence properties
\begin{eqnarray*}
\lim_{t\to\infty}\sup_{x\in\Ri^d}\int_{\Ri^d}dy\,
|K_t(x\,;y)-{\widehat K}_t(x\,;y)|&=&0\\[5pt]
\lim_{t\to\infty}\sup_{y\in\Ri^d}\int_{\Ri^d}dx\,
|K_t(x\,;y)-{\widehat K}_t(x\,;y)|&=&0
\end{eqnarray*}
of the kernel.

\section{Spectral refinement}

The homogenization process described in Section~2 shows that the rescaled
versions
$H^{(m),z}$ of the $z$-periodic operators converge in a very strong sense.
Hence one has good estimates on their spectral properties especially if the
operators
are self-adjoint.
The basis of the discussion was the observation that the operators are all
periodic and
the scaling decreased the period.
But this means that the decomposition theory of Section~2 has further
refinements and
the limit of the spectrum is related to this refinement.
In this section we analyze the details of the refinement process in a
slightly more
general context.

We now consider self-adjoint operators, i.e., we assume $C=(c_{ij})$ is
hermitian,
$c_i'=-{\overline c}_i$ and $c_0$ is real. 
This ensures that $H$ is self-adjoint on $L_2(\Ri^d)$ and the operators
$H_z$ with
$z$-periodic boundary conditions are self-adjoint on $L_2(\Ii^d)$.
We further assume that $H$ has a H\"older continuous Gaussian kernel so all the
preceding results apply.
Let $H^M$ and $H_z^M$ be the operators obtained from $H$ and $H_z$ by
the replacements
\[
c_{ij}(x)\to c_{ij}^M(x)=c_{ij}(Mx)
\]
and $c_i(x)\to c_i^M(x)=c_i(Mx)$ etc.\ where $M$ is a $d\times d$-matrix
with integer
coefficients and $N=|\det M|>1$. Subsequently, in the case that $M$ is a
multiple of the
identity, we consider sequences of operators obtained by iterating this map.

The new coefficients still satisfy the periodicity condition
$c^M_{ij}(x+n)=c^M_{ij}(x)$ for all $n\in\Zi^d$ but they also satisfy the
further
periodicity property
\begin{equation}
c^M_{ij}(x+M^{-1}n)=c^M_{ij}(x)\label{eper4.1}
\end{equation}
for  $M^{-1}n\in M^{-1}(\Zi^d)$ and this leads to an immediate refinement
of the Zak
decomposition.
Define a generalized Zak transform $Z_M\colon L_2(\Ri^d)\mapsto L_2(\Ti^d\times
M^{-1}(\Ii^d))$ by
\begin{equation}
(Z_Mf)(z,u)=\sum_{n\in\Zi^d}z^nf(u-M^{-1}n)\;\;\;.\label{eper4.2}
\end{equation}
One calculates as before that $Z_M$ is a unitary operator when
$M^{-1}(\Ii^d)$ is
equipped with Lebesgue measure inherited from $\Ri^d$.
One may again extend the definition of $Z_Mf$ from $u\in M^{-1}(\Ii^d)$ to
general
$u\in\Ri^d$ by periodicity and one then has the identity 
\begin{equation}
(Z_Mf)(z,x+M^{-1}n)=z^n(Z_Mf)(z,x)\label{eper4.3}
\end{equation}
for all $z\in\Ti^d$, $x\in\Ri^d$ and $n\in\Zi^d$.
This leads, as before, to a decomposition
\begin{equation}
L_2(\Ri^d)=(2\pi)^{-d}\int^\oplus_{\Ti^d}|dz|\,{\cal H}_M(z)\label{eper4.4}
\end{equation}
where ${\cal H}_M(z)\simeq L_2(M^{-1}(\Ii^d))$.

If $A$ is an operator on $L_2(\Ri^d)$ which commutes with the action of
$M^{-1}\Zi^d$ then $A$ or, more precisely, $Z_MAZ_M^*$ has a decomposition
\[
A=(2\pi)^{-d}\int^\oplus_{\Ti^d}|dz|\,A_{M,z}\;\;\;.
\]
 $\Zi^d\subseteq M^{-1}\Zi^d$ the operator $A$ also commutes with the action of
$\Zi^d$ and has the original Zak decomposition
\[
A=(2\pi)^{-d}\int^\oplus_{\Ti^d}|dz|\,A_z\;\;\;.
\]
In particular if $A=H^M$ then $H^M$ has both decompositions into $H^M_{M,z}$
and into $H^M_z$.
We now study the connection between the two decompositions and we first
demonstrate
 that the $Z_M$-decomposition is a refinement of the $Z$-decomposition,
i.e., each of the Hilbert spaces ${\cal H}(z)$ which occurs in the
$Z$-decomposition
of $L_2(\Ri^d)$ decomposes into a direct sum of $N$ of the Hilbert spaces ${\cal
H}_M(z)$ in such a way that the corresponding components of $H^M$, and the
semigroup
$S^M$ generated by $H^M$, decompose accordingly.
Adopting the notation
\begin{equation}
z^M=(e^{i\theta_1},\ldots,e^{i\theta_d})^M
=(e^{i\sum_{j=1}^dM_{j1}\theta_j},\ldots,e^{i\sum_{j=1}^dM_{jd}\theta_j})
\label{eper4.5}
\end{equation}
this decomposition can be described as follows:
\begin{prop}\label{pper4.1}
Each Hilbert space ${\cal H}(z)\,(\simeq L_2(\Ii^d))$ in the Zak
decomposition of
$L_2(\Ri^d)$ has the decomposition
\begin{equation}
{\cal H}(z)=\bigoplus_{{w\in\Ti^d}\atop{w^M=z}}
{\cal H}_M(w)\;\;\;.\label{eper4.6}
\end{equation}
More specifically if $F=Zf$ and $F_M=Z_Mf$ with $f\in L_2(\Ri^d)$ then 
\begin{equation}
F(z\,,x)=N^{-1}\sum_{{w\in\Ti^d}\atop{w^M=z}}F_M(w\,,x)\label{eper4.7}
\end{equation}
for all  $z\in\Ti^d$ and $x\in\Ri^d$.

Conversely, ${\cal H}_M(z)$ is embedded into ${\cal H}(z^M)$ by
\begin{equation}
F_M(z\,,x)=\sum_{p\in\Zi^d/M\Zi^d}z^{-p}
F(z^M\,,x+M^{-1}p)\label{eper4.8}
\end{equation}
for all  $z\in\Ti^d$ and $x\in\Ri^d$.
\end{prop}
\begin{remarkn}\label{rper4.1} Formula (\ref{eper4.7}) states that if
$u\in\Ii^d$
and $i(u)$ is the unique element of $\Zi^d$  such that 
$u-M^{-1}(i(u))\in M^{-1}(\Ii^d)$
then
\begin{equation}
F(z\,,u)=N^{-1}\sum_{{w\in\Ti^d}\atop{w^M=z}}
w^{i(u)}F_M(w\,,u-i(u))\label{eper4.9}
\end{equation}
which implies (\ref{eper4.6}).

\end{remarkn}

\begin{remarkn}\label{rper4.2}If $w^M=z$ the projection $P$ from ${\cal
H}(z)$ onto the
subspace ${\cal H}_M(w)$ may be written explicitly as 
\begin{equation}
(PF)(u)=N^{-1}\sum_{p\in\Zi^d/M\Zi^d}w^{-p}
F(z\,,u+M^{-1}p)\;\;\;.\label{eper4.10}
\end{equation}
(See Remark~\ref{rper4.3} below.)
\end{remarkn}
\bigskip

\noindent{\bf Proof of Proposition~\ref{pper4.1}} By definition
\begin{eqnarray}
F(z\,,u)=(ZZ_M^{-1}F_M)(z\,,u)
&=&\sum_{n\in\Zi^d}z^n(Z_M^{-1}F_M)(u-n)\nonumber\\[5pt]
&=&\sum_{n\in\Zi^d}z^n(2\pi)^{-d}
\int_{\Ti^d}|dw|\,F_M(w\,,u-n)\nonumber\\[5pt]
&=&\sum_{n\in\Zi^d}z^n(2\pi)^{-d}
\int_{\Ti^d}|dw|\,w^{-Mn}F_M(w\,,u)\;\;\;.
\label{eper4.11}
\end{eqnarray}
Next define the mean $R$ by
\begin{equation}
(RF_M)(w\,,u)=N^{-1}\sum_{{v\in\Ti^d}\atop{v^M=w}}
F_M(v\,,u)\label{eper4.12}
\end{equation}
then $R$ is the adjoint of the operator $f(z)\mapsto f(z^M)$ and hence
\begin{equation}
F(z\,,u)=\sum_{n\in\Zi^d}z^n(2\pi)^{-d}
\int_{\Ti^d}|dw|\,w^{-n}(RF_M)(w\,,u)
=(RF_M)(z\,,u)\label{eper4.13}
\end{equation}
and (\ref{eper4.7}) is established.
The decomposition (\ref{eper4.6}) follows automatically (see
Remark~\ref{rper4.1}).

The converse formula (\ref{eper4.8}) is a consequence of the calculation
\begin{eqnarray}
F_M(z\,,u)&=&(Z_MZ^{-1}F)(z\,,u)\nonumber\\[5pt]
&=&\sum_{n\in\Zi^d}z^n(Z^{-1}F)(u-M^{-1}n)\nonumber\\[5pt]
&=&\sum_{n\in\Zi^d}z^n(2\pi)^{-d}
\int_{\Ti^d}|d\zeta|\,F(\zeta\,,u-M^{-1}n)
\nonumber\\[5pt]
&=&\sum_{k\in\Zi^d}\sum_{y\in M^{-1}\Zi^d/\Zi^d}
z^{(Mk+My)}(2\pi)^{-d}\int_{\Ti^d}|d\zeta|\,F(\zeta\,,u-y-k)
\nonumber\\[5pt]
&=&\sum_{k\in\Zi^d}\sum_{y\in M^{-1}\Zi^d/\Zi^d}
z^{Mk}z^{My}(2\pi)^{-d}
\int_{\Ti^d}|d\zeta|\,\zeta^{-k}F(\zeta\,,u-y)
\nonumber\\[5pt]
&=&\sum_{y\in M^{-1}\Zi^d/\Zi^d}
z^{My}F(z^M\,,u-y)
\nonumber\\[5pt]
&=&\sum_{p\in \Zi^d/M\Zi^d}
z^{-p}F(z^M\,,u+M^{-1}p)
\label{eper4.14}
\end{eqnarray}
and this completes the proof.\hfill\Box
\bigskip
\begin{cor}\label{cper4.1}
The $z$-component $H^M_z$ of $H^M$ in the $Z$-decomposition has the
further decomposition
\begin{equation}
H^M_z=\bigoplus_{{w\in\Ti^d}\atop{w^M=z}}H^M_{M,w}
\label{eper4.15}
\end{equation}
corresponding to the decomposition {\rm (\ref{eper4.6})} of ${\cal H}(z)$
so that the
spectrum over the point $z$ in the
$Z$-decomposition is the union of the spectra of the
$H^M_{M,w}$ over the $w$ with $w=z^M$ in the $Z_M$-decomposition.
\end{cor}
This follows from Proposition~\ref{pper4.1} and the fact that $H^M$ is
decomposable with
respect to both Zak decompositions.

\begin{remarkn}\label{rper4.3}The operator decomposition (\ref{eper4.15})
refers to the
Hilbert space decomposition (\ref{eper4.6}) which in turn can be made more
precise by
the identification of the corresponding isometric embeddings 
$S_w\colon F_M(w\,,\cdot)\in{\cal H}_M(w)\mapsto 
N^{-1/2}F_M(w\,,\cdot)\in{\cal H}(z)$
where $z=w^M$.
The operators $S_w$ are well defined since
\[
F_M(w\,,x+n)=F_M(w\,,x+M^{-1}Mn)=w^{Mn}F_M(w\,,x)=
z^nF_M(w\,,x)
\]
for all $x\in\Ri^d$ and $n\in\Zi^d$.
The adjoint operator
$S^*_w\colon {\cal H}(z)\mapsto {\cal H}_M(w)$ is given by
\[
(S^*_wF)(z)=N^{-1/2}\sum_{p\in\Zi^d/M\Zi^d}w^{-p}
F(z^M\,,x+M^{-1}p)
\;\;\;.
\]
It can be checked that the projections $P_w=S_wS^*_w$ on ${\cal H}(z)$ are
mutually
orthogonal and
\[
\sum_{w:w^M=z}P_w=I_{{\cal H}(z)}
\;\;\;.
\]
The formula for $P_w$ is the finite Fourier transform in (\ref{eper4.10}).
We may use this and (\ref{eper4.15}) to identify the eigenvectors
corresponding to the
refined spectrum.
Let $\lambda(z)$ be in the spectrum of $ZH^MZ^{-1}$ with eigenvector
$F\in{\cal H}(z)$.
The projections $P_w=S_wS^*_w$ commute with $ZH^MZ^{-1}$.
The $w\in\Ti^d$ such that $F\in(Z_MH^MZ_M^{-1})(w)$ are exactly those $w$
(generically
only one) such that $S^*_wF\neq0$.
This happens if and only if the corresponding conjugacy class in $\Zi^d/M\Zi^d$
contributes to the Fourier series expansion of $F$.
\end{remarkn}

In order to analyze the spectrum in more detail we specialize the
transformation $M$. 
We assume $M=N I$ where $N\in\{2,3,\ldots\}$ so that $\det M=N^d$ and we
use $N$ as
index instead of $M$, i.e., we write $Z_N$ instead of $Z_M$ etc.
This special choice of $M$ corresponds to a rescaling of the type considered in
Section~3 and hence the semigroups $S^N$ generated by the rescaled
operators $H^N$
converge to the semigroup $\widehat S$ generated by the  homogenization
$\widehat H$ of
$H$. 
We expect that for general transformations $M$ of the foregoing type the
corresponding semigroups $S^M$ will have good convergence properties if all the
eigenvalues of $M$ have modulus strictly larger than one but it appears
difficult to
identify their limit. Hence we restrict attention to the special case
$M=NI$ and discuss
the links between the spectral refinement and the
homogenization limit which we obtained in Section~3.

The special choice of $M$ now gives 
\begin{equation}
(Z_Nf)(z\,,u)=\sum_{n\in\Zi^d}z^nf(u-n/N)\label{eper4.16}
\end{equation}
for $u\in N^{-1}\Ii^d$ and (\ref{eper4.9}) takes the form
\begin{equation}
F(z\,,u+N^{-1}i)=N^{-d}\sum_{{w\in\Ti^d}\atop{w^N=z}}w^i 
F_N(w\,,u)\label{eper4.17}
\end{equation}
for $u\in N^{-1}\Ii^d$, $i=(i_1,\ldots,i_d)$ with $i_k\in\{0,1,\ldots, N-1\}$.

Let $H^N_z$ and $H^N_{N,z}$ denote the components of the rescaled operator
$H^N$ in
the $Z$ and $Z_N$ decompositions, respectively.
For all $z\in\Ti^d$ let $\lambda_n(z)$ denote the eigenvalues of $H_z$ on
$L_2(\Ii^d)$
arranged in increasing order repeated according to multiplicity.
Now consider pure second-order
operators, $H=-\sum^d_{i,j=1}\partial_ic_{ij}\partial_j$.
Then by scaling the eigenvalues of $H^N_{N,z}$ on $L_2(N^{-1}\Ii^d)$ are
$N^2\lambda_n(z)$. 
As a consequence of Corollary~\ref{cper4.1} one
draws the following  conclusion about the spectrum.
\begin{cor}\label{cper4.2}
Let $H$ be a pure second-order, self-adjoint, elliptic operator and  let
$\lambda_n(z)$
denote the eigenvalues of $H_z$ on $L_2(\Ii^d)$. 
The eigenvalues, counted with
multiplicity, of the rescaled operator
\[
H^N_z=-\sum^d_{i,j=1}\partial_ic_{ij}(N\cdot)\partial_j
\]
with $z$-periodic boundary conditions on $L_2(\Ii^d)$ form the set of
$N^2\lambda_n(w)$
with $n=0,1,\ldots$ and $w\in\Ti^d$ with $w^N=z$.
\end{cor}

If we now take the limit $N\to\infty$ the semigroups $S^{N,z}$, generated
by $H^N_z$ on
$L_2(\Ii^d)$, converge in the trace norm, by Theorem~\ref{tper3.1},
uniformly for $t$ in
compact intervals of $\langle0,\infty\rangle$ to the semigroup ${\widehat S}_z$
generated by the homogenization ${\widehat H}_z$ of $H_z$. 
Therefore the eigenvalues
converge pointwise to those of
${\widehat H}_z$ and the finite-dimensional eigenprojections converge in norm.
But the spectrum of ${\widehat H}_z$ with
$z=(e^{i\theta_1},\ldots,e^{i\theta_d})$ is
easily computed.
For $z$ fixed set $\varphi_n(u)=e^{i(\theta-n).u}$ for $n\in\Zi^d$ then
\[
{\widehat H}_z\varphi_n=\sum^d_{i,j=1}{\widehat
c}_{ij}(n_i-\theta_i)(n_j-\theta_j)\varphi_n
\]
so the eigenvalues of ${\widehat H}_z$, counted with multiplicity, are
\begin{equation}
\{\langle (n-\theta),{\widehat C}(n-\theta)\rangle\,
:\,n\in\Zi^d\,\}\label{eper4.37}
\;\;\;.
\end{equation}
Therefore one has the following asymptotic identification of the spectrum
of $H_z$.
\begin{cor}\label{cper4.3}
If $\lambda_n(z)$ denotes the eigenvalues of $H_z$ then
\[
\lim_{N\to\infty}\{N^2\lambda_n(w):\,w^N=z\,,\,n=0,1,\ldots\,\}=
\{\langle (n-\theta),{\widehat C}(n-\theta)\rangle
\,:\,n\in\Zi^d\,\}
\]
where the limit is  in the sense of pointwise convergence of the ordered sets.
\end{cor}

The rate of convergence of the eigenvalues can be estimated further by the
trace norm
estimate (\ref{eper3.4}).

\section{Concluding remarks}

Adopt the general assumptions in the beginning of Section~4, so that $H$ is
self-adjoint on
$L_2(\Ri^d)$.
Again let $\lambda_n(z)$, $z\in\Ti^d$, denote the eigenvalues of $H_z$ on
$L_2(\Ii^d)$ arranged in
increasing order repeated according to multiplicity.
It then follows from Corollary~\ref{cper2.1} and the minimax principle that
the functions $z\mapsto
\lambda_n(z)$ are continuous on $\Ti^d$ and can even be shown to be
piecewise analytic
in a suitable sense \cite{Kat1}, Theorem~VII.7.18, \cite{Don}, Lemma~2.1.
It then follows, see for example \cite{RS4}, Theorem~XIII.85, that the
spectrum of $H$ on
$L_2(\Ri^d)$ is the union of the closed intervals
\[
B_n=\{\lambda_n(z):\,z\in\Ti^d\,\}
\;\;\;.
\]
These are called the bands of the spectrum and they may, or may not,
overlap and in
principle they could, or could not, contain points which are eigenvalues
(corresponding to
open sets in
$\Ti^d$ where $\lambda_n$ is constant).
Because of piecewise analyticity the spectrum of $H$, apart from the
possible eigenvalues, is
absolutely continuous.
Corollary~\ref{cper4.3} shows among other things that possible gaps in the
bands of
$H^N$ disappear as $N\to\infty$ and the spectrum converges to the wellknown
spectrum
$[0,\infty\rangle$ of ${\widehat H}$ with multiplicities given by Fourier
analysis.

One instance where $H$ has no eigenvalues is  the special case  of
Schr\"odinger operators,
\begin{equation}
H=-\Delta +V \;\;\;, \label{eper5.1}
\end{equation}
with periodic potentials $V$.
There has been a great deal of interest in these operators since the
seminal work of Bloch
\cite{Blo}  (see, for example,  \cite{DT} \cite{Eas} \cite {RS4} \cite{Skr}
\cite{Sun} \cite{Wilc}).
It is known that the spectrum of $H$ on $L_2(\Ri^d)$ is absolutely
continuous (see, for
example, \cite{Thom} \cite{Wilc} and \cite{RS4}, Theorem~XIII.100).
The argument is roughly the following.
If $\lambda>0$ and $t>0$ are fixed, the operators $(\lambda I+H_z)^{-1}$
and $S^z_t$ have analytic
continuations to $0<|z|<+\infty$ by the results in Section~2.
Introducing $k\in\Ci^d$ by $z=\exp(ik)$, these functions are entire
analytic in $k$ with values in
the trace class operators, by Corollary~\ref{cper2.1}.
Now if the $\lambda_n$, viewed as functions of $k$, are constant in an open
set then $H_z$ has a
constant eigenvalue for $z\in\Ci^d\backslash\{0\}$ by \cite{Kat1} \cite{Don}.
But then the norms $\|(\lambda I+H_z)^{-1}\|_{L_2(\Ii^d)}$ and
$\|S^z_t\|_{L_2(\Ii^d)}$ would have
positive lower bounds uniformly in  $z\in\Ci^d\backslash\{0\}$.
If, however, $H$ is a Schr\"odinger operator then Thomas \cite{Thom} shows that 
$\|(\lambda I+H_z)^{-1}\|_{L_2(\Ii^d)}\to 0$
as $k\to \infty$ through a line  in the imaginary direction.
The estimates in Section~2 could possibly be elaborated to reach a similar
conclusion  for
general self-adjoint elliptic operators $H$.

Returning to the Schr\"odinger operator (\ref{eper5.1}),
we have seen that its spectrum is the union of bands 
\[
B_n(V)=\{\lambda_n(z):\,z\in\Ti^d\,\}\;\;\;.
\]
If $d=1$ the $n$-th band is located in a neighbourhood of $\pi n^2$ and the
bands may or may not
overlap depending on the nature of $V$.
In general the length of the gaps between bands, if such a gap exists,
tends to zero as $n\to
\infty$ and if ${\overline V}=\int^1_0dx\,V(x)$ is the mean of $V$ then the
length of each
gap is dominated by $2(\int^1_0dx\,(V(x)-{\overline V})^2)^{1/2}$ in the
region $\langle
{\overline V},\infty\rangle$.
Thus  there  are no gaps  if $V$ is constant as we can also see from the
direct computation 
 $\{\lambda_n(e^{i\theta}):\,n\in\Ni\,\}=
\{(k+\theta)^2+V:\,k\in\Zi\,\}$.
Conversely it can be proved that if there are no gaps then $V$ is constant.
There are examples where all gaps are open.
If there is precisely one gap then $V$ is a Weierstrass elliptic function.
If there are only finitely many gaps then $V$ is real analytic and if all
the odd gaps are absent
then $V$ has period $1/2$ (see \cite{Eas} and \cite{RS4}, Section~VIII.16,
for a discussion of
these results).

If $d=2$ it is established in \cite{DT} that there exist positive $c$ and
$C$ such that
\[
[4\pi n-cn^{1/4},4\pi n+cn^{1/4}]\subseteq B_n(V)\subseteq
[4\pi n-Cn^{1/3},4\pi n+Cn^{1/3}]
\;\;\;.
\]
Hence the spectrum may have a finite number of open gaps near zero but the
spectrum always
contains all sufficiently large $\lambda\in\Ri$.
If $d= 3$ this qualitative picture remains.
If one defines a multiplicity function $M$ by
\[
M(\lambda)=\#\{n:\,\lambda\in B_n(V)\,\}
\]
it is proved in \cite{Skr} that there exist positive $c$ and $C$ such that
\[
c\,\lambda^{1/2}< M(\lambda) < C\,\lambda^{3/4}
\]
for large $\lambda$ and it is conjectured that $\lambda^{3/4}$ may be
replaced by
$\lambda^{1/2+\varepsilon}$ in the upper estimate.
In this case it can also be shown that if $\|V\|_\infty< D$ where $D$ is a
universal constant
there are no gaps in the spectrum.

If $V$ is replaced by $V_\varepsilon=V(\varepsilon^{-1}\cdot)$ it follows from
Corollary~\ref{cper3.1} that $H_\varepsilon=-\Delta+V_\varepsilon$
converges to ${\widehat
H}=-\Delta+{\overline V}$ in the norm resolvent sense
with 
\begin{equation}
{\overline V}=\int_{\Ii^d}dx\,V(x)\label{eper5.2}
\;\;\;.
\end{equation}
Since ${\overline V}$ is constant the spectral properties of ${\widehat
H}$ are the same as those of $-\Delta$ and we may formulate a corollary
similar to
Corollary~\ref{cper4.3} although the eigenvalues of $H_{N,z}$ can no longer
be computed from
those of $H_z$ by simple scaling properties.
Let us finally remark that if $H_\varepsilon$ is an operator of the form
\begin{equation}
H_\varepsilon=H_0+V_\varepsilon\label{eper5.3}
\end{equation}
where $H_0$ is any elliptic operator coming from a form of the type
(\ref{per2.1}), with no
periodicity required, with a kernel satisfying the bounds (\ref{per2.13})
and (\ref{per2.14}) then
$H_\varepsilon\to {\widehat H}=H_0+{\overline V}$ in the norm resolvent sense.
This can be established by the following direct argument.

The Duhamel expansion
\[
S^\varepsilon_t-{\widehat S}_t=\int^t_0ds\,S^\varepsilon_s
({\overline V}-V_\varepsilon){\widehat S}_{t-s}
\]
for the corresponding semigroups is valid on $L_2(\Ri^d)$ and for fixed $t$
the kernel of the
operator in the integrand is 
\begin{eqnarray*}
K(s\,:x\,;y)&=&\int_{\Ri^d}du\,K^\varepsilon_s(x\,;u)
(({\overline V}-V(\varepsilon^{-1}u))
{\widehat K}_{t-s}(u\,;y)\\[5pt]
&=&\sum_{n\in\Zi^d}\varepsilon^d
\int_{n+\Ii^d}du\,K^\varepsilon_s(x\,;\varepsilon u)
(({\overline V}-V(u)) {\widehat K}_{t-s}(\varepsilon u\,;y)
\;\;
\end{eqnarray*}
Now we use the H\"older continuity (\ref{per2.14}) to estimate each of the
terms in the sum by
$O(\varepsilon^\nu)$ when $\varepsilon u$ is close to $x$ or $y$  and use
the Gaussian bounds
(\ref{per2.13}) to estimate in terms of $\|{\overline V}-V\|_\infty$ times
a Gaussian when
$\varepsilon u$ is far from $x$ and $y$.
These estimates are singular in $s$ when it is near $0$ or $t$ but this can
be remedied by first
estimating the Duhamel integral near these limits in terms of $\|
V-V_\varepsilon\|_\infty$.
These estimates lead to the conclusion that 
\[
\lim_{\varepsilon\to0}\sup_{x\in\Ri^d}\int_{\Ri^d}dy\,|K(s\,:x\,;y)|=0
\]
uniformly for $s$ in compact intervals of $\langle0,t\rangle$.
One then concludes that one has norm convergence of the semigroups,
$\|S^\varepsilon_t-{\widehat
S}_t\|_{L_p(\Ri^d)\to L_p(\Ri^d)}\to0$ as
$\varepsilon\to 0$ for each $p\in[1,\infty]$, uniformly for $t$ in compact
intervals of
$\langle0,\infty\rangle$.

\section*{Acknowledgements}

This work was carried out whilst the first two authors were visiting the
Centre for 
Mathematics and its Applications  at the ANU with financial support of the
Centre.
The first author also obtained travel support from the Norwegian Research
Council.
All three authors are indebted to Tom ter Elst for a critical reading of
the manuscript and
several helpful suggestions.

\bibliography{refbib}

\begin{thebibliography}{ZKON79}

\bibitem[Ati76]{Ati}
{\sc Atiyah, M.~F.}, Elliptic operators, discrete groups and von Neumann
  algebras.
\newblock {\em Ast{\'e}risque} {\bf 32-33} (1976).

\bibitem[Aus96]{Aus1}
{\sc Auscher, P.}, Regularity theorems and heat kernel for elliptic operators.
\newblock {\em J. London Math.\ Soc.} {\bf 54} (1996),  284--296.

\bibitem[ACT96]{ACT}
{\sc Auscher, P., Coulhon, T., {\rm and} Tchamitchian, P.}, Absence de principe
  du maximum pour certaines \'equations paraboliques complexes.
\newblock {\em Colloq.\ Math.} {\bf 71} (1996),  87--95.

\bibitem[AMT94]{AMT1}
{\sc Auscher, P., McIntosh, A., {\rm and} Tchamitchian, P.}, Heat kernels of
  second order complex elliptic operators and their applications.
\newblock Research Report 94---164, Maquarie University, Sydney, Australia,
  1994.

\bibitem[BBJR95]{BBJR}
{\sc Batty, C. J.~K., Bratteli, O., J{\o}rgensen, P. E.~T., {\rm and} Robinson,
  D.~W.}, Asymptotics of periodic subelliptic operators.
\newblock {\em J. Geom.\ Anal.} {\bf 5} (1995),  427--443.

\bibitem[BLP78]{BLP}
{\sc Bensoussan, A., Lions, J.~L., {\rm and} Papanicolaou, G.}, {\em Asymptotic
  analysis for periodic structures}, vol.\ 5 of Studies in Mathematics and its
  Applications.
\newblock North-Holland, Amsterdam etc., 1978.

\bibitem[Blo28]{Blo}
{\sc Bloch, F.}, Uber die quantenmechanik der elektronen in kristallgittern.
\newblock {\em Zeit.\ f{\"u}r Physik} {\bf 52} (1928),  555--600.


\bibitem[BR87]{BR1}
{\sc Bratteli, O., {\rm and} Robinson, D.~W.}, {\em Operator algebras and
  quantum statistical mechanics}, vol.\ 1.
\newblock Second edition. Springer-Verlag, New York etc., 1987.

\bibitem[DT82]{DT}
{\sc Dahlberg, B. E.~J., {\rm and} Trubowitz, E.}, A remark on two dimensional
  periodic potentials.
\newblock {\em Comm.\ Math.\ Helvetici} {\bf 57} (1982),  130--134.

\bibitem[Dau92]{Dau}
{\sc Daubechies, I.}, {\em Ten lectures on wavelets}.
\newblock CBMS Regional Conference Series in Mathematics 61. Amer.\ Math.\
  Soc., Providence, 1992.

\bibitem[Dix69]{Dix}
{\sc Dixmier, J.}, {\em les alg\`ebres d'op\'erateurs dans l'espace
  Hilbertien}.
\newblock Cahiers scientifiques 25. Gauthier-Villars, Paris, 1969.

\bibitem[Don81]{Don}
{\sc Donnelly, H.}, On $L_2$-Betti numbers for abelian groups.
\newblock {\em Can.\ Math.\ Bull.} {\bf 24} (1981),  91--95.


\bibitem[Eas73]{Eas}
{\sc Eastham, M. S.~P.}, {\em The spectral theory of periodic differential
  equations}.
\newblock Scottish Academic Press, Edinburgh and London, 1973.

\bibitem[ER96]{ER15}
{\sc Elst, A. F.~M. ter, {\rm and} Robinson, D.~W.}, Second-order subelliptic
  operators on Lie groups I: complex uniformly continuous principal
  coefficients.
\newblock Research Report MRR 035-96, The Australian National University,
  Canberra, Australia, 1996.

\bibitem[ER97]{ER20}
\leavevmode\vrule height 2pt depth -1.6pt width 23pt, High order
  divergence-form elliptic operators on Lie groups.
\newblock {\em Bull.\ Austral.\ Math.\ Soc.} {\bf 55} (1997),  335--348.

\bibitem[FS86]{FaS}
{\sc Fabes, E.~B., {\rm and} Stroock, D.~W.}, A new proof of Moser's parabolic
  Harnack inequality using the old ideas of Nash.
\newblock {\em Arch.\ Rat.\ Mech.\ and Anal.} {\bf 96} (1986),  327--338.

\bibitem[Gia83]{Gia1}
{\sc Giaquinta, M.}, {\em Multiple integrals in the calculus of variations and
  nonlinear elliptic systems}.
\newblock Annals of Mathematics Studies 105. Princeton University Press,
  Princeton, 1983.

\bibitem[Gia93]{Gia}
\leavevmode\vrule height 2pt depth -1.6pt width 23pt, {\em Introduction to
  regularity theory for nonlinear elliptic systems}.
\newblock Lectures in Mathematics ETH Z{\"u}rich. Birkh{\"a}user Verlag, Basel
  etc., 1993.

\bibitem[Gio57]{DG}
{\sc Giorgi, E.~D.}, Sulla differenziabilit\`a e l'analiticit\`a delle
  estremali degli integrali multipli regolari.
\newblock {\em Mem.\ Accad.\ Sci.\ Torino cl.\ Sci.\ Fis.\ Mat.\ Nat.} {\bf 3}
  (1957),  25--43.

\bibitem[Kat61]{Kat2}
{\sc Kato, T.}, Fractional powers of dissipative operators.
\newblock {\em J. Math.\ Soc.\ Japan} {\bf 13} (1961),  246--274.

\bibitem[Kat84]{Kat1}
\leavevmode\vrule height 2pt depth -1.6pt width 23pt, {\em Perturbation theory
  for linear operators}.
\newblock Second edition, Grundlehren der mathematischen Wissenschaften 132.
  Springer-Verlag, Berlin etc., 1984.

\bibitem[Nas58]{Nash}
{\sc Nash, J.}, Continuity of solutions of parabolic and elliptic equations.
\newblock {\em Amer.\ J. Math.} {\bf 80} (1958),  931--954.

\bibitem[Pow67]{Pow}
{\sc Powers, R.~T.}, Representations of uniformly hyperfinite algebras and
  their associated von Neumann rings.
\newblock {\em Ann. Math.} {\bf 86} (1967),  138--171.

\bibitem[RS78]{RS4}
{\sc Reed, M., {\rm and} Simon, B.}, {\em Methods of modern mathematical
  physics IV. Analysis of operators}.
\newblock Academic Press, New York etc., 1978.

\bibitem[Skr85]{Skr}
{\sc Skriganov, M.~M.}, The spectrum band structure of the three-dimensional
  Schr{\"o}dinger operator with periodic potential.
\newblock {\em Invent.\ Math.} {\bf 80} (1985),  107--121.

\bibitem[Sun88]{Sun}
{\sc Sunada, T.}, Fundamental groups and Laplacians.
\newblock In  {\em Geometry and analysis on manifolds}, 
  Lecture Notes in Mathematics 1339.
  Springer-Verlag, Berlin etc., 1988.

\bibitem[Tho73]{Thom}
{\sc Thomas, L.~E.}, Time dependent approach to scattering from impurities in a
  crystal.
\newblock {\em Commun.\ Math.\ Phys.} {\bf 33} (1973),  335--343.

\bibitem[Wil78]{Wilc}
{\sc Wilcox, C.~H.}, Theory of Bloch waves.
\newblock {\em J. d'Analyse Math.} {\bf 33} (1978),  146--167.

\bibitem[ZKON79]{ZKON}
{\sc Zhikov, V.~V., Kozlov, S.~M., Oleinik, O.~A., {\rm and} Ngoan, K.~T.},
  Averaging and $G$-convergence of differential operators.
\newblock {\em Russian Math.\ Surveys} {\bf 34} (1979),  569--147.

\end{thebibliography}

\end{document}